\documentclass[11pt,thmsa]{article}
\usepackage{amsmath}
\usepackage{amsfonts}
\usepackage{amssymb}
\usepackage{graphicx}
\newtheorem{theorem}{Theorem}[section]

\newtheorem{problem}[theorem]{Problem}

\newtheorem{definition}[theorem]{Definition}

\newtheorem{lemma}[theorem]{Lemma}

\begin{document}
\title{Recent progress on homogeneous Finsler spaces with positive curvature\thanks{Supported by NSFC (no. 11671212, 51535008), SRFDP of China}}
\author{Shaoqiang Deng$^1$, Ming Xu$^2$\thanks{Corresponding author} \\
\\
$^1$School of Mathematical Sciences and LPMC\\
Nankai University\\
Tianjin 300071, P. R. China\\
E-mail: dengsq@nankai.edu.cn\\
\\
 $^2$College of Mathematics\\
Tianjin Normal University\\
 Tianjin 300387, P. R. China\\
 Email: mgmgmgxu@163.com.}

\date{}
\maketitle

\begin{abstract}
This is a survey paper on our recent works concerning  the classification of
positively curved homogeneous Finsler spaces, and some related
topics. At the final part, we present some open problems in this field.


\end{abstract}

\section{Introduction}
The purpose of this article is to survey some recent works on the classification of homogeneous Finsler spaces with positive flag curvature, which is the most important topic in the more generalized subject related to the classification of Finsler spaces with positive flag curvature. The same problem in Riemannian geometry has been one of the central problems in differential geometry. In the homogeneous case, a complete classification was achieved by the works of Berger, Wallach, Aloff-Wallach, B\'{e}rard Bergery, Wilking, Xu-Wolf and Wilking-Ziller; See \cite{AW75, BB76, Ber61, Wallach1972, Wi1999, WZ2015, XW2015}.

The study of the above problem in the Finsler case was initiated  by S. Deng and Z. Hu in \cite{HD11}, where they classified homogeneous Randers metrics with positive flag curvature and vanishing S-curvature. Note that their classification is also valid for homogeneous $(\alpha,\beta)$-spaces with positive flag curvature and vanishing S-curvature \cite{XD2015}.
Recently,  big progress has been made on the classification with more generality.
In \cite{XD-Normal-homogeneous-Finsler-spaces},  the authors of this paper  classified positively curved
normal homogeneous Finsler spaces, generalizing the classical results of \cite{Ber61}. In the joint work of the authors with L. Huang and Z. Hu \cite{XDHH2014}, we classified even-dimensional positively curved homogeneous Finsler spaces, generalizing the results of \cite{Wallach1972}. In the odd-dimensional case, big steps have also been made towards a  classification
by Xu-Deng \cite{XD2016} and \cite{XZi2016}. However, a complete classification
has  not been achieved yet, due to some technical difficulties.

In Section 2, we give some preliminaries in Finsler geometry. In Section 3, we state the problems related to the topics of this paper. Section 4 is devoted to
presenting the techniques which we have used in our study. In Sections 5, 6, and 7, we recall the main progress in the study of positively curved Finsler spaces. In Section 8, we recall some results on negatively curved Finsler spaces. Finally, in Section 9, we pose some open problems related to the topics of this paper.

\section{Preliminaries}

In this section, we collect some preliminaries of Finsler geometry which will be used in this paper. See \cite{BCS2000}\cite{CS2005}\cite{Sh2001} for more details.

\subsection{Finsler metrics and Minkowski norms}

A {\it Finsler metric} on a smooth manifold $M$ is a continuous function
$F:TM\rightarrow [0,+\infty)$ satisfying the following conditions:
\begin{description}
\item{\rm (1)} Regularity: the restriction of $F$ to the slit tangent bundle
$TM\backslash 0$ is a positive smooth function.
\item{\rm (2)} Homogeneity: for any $x\in M$, $y\in T_xM$ and $\lambda\geq 0$, $F(x,\lambda y)=\lambda F(x,y)$.
\item{\rm (3)} Strong Convexity: given any standard local coordinates on $TM$, $x=(x^i)\in M$
and $y=y^j\partial_{x^j}\in T_xM$, the Hessian matrix
$$(g_{ij}(x,y))=\left(\frac12\frac{\partial^2}{\partial y^i\partial y^j}[F^2(x,y)]\right)$$
is positive definite when $y\neq 0$.
\end{description}
We will call $(M,F)$ a {\it Finsler manifold} or a {\it Finsler space}. For any nonzero tangent vector $y\in T_xM$, The Hessian matrix $(g_{ij}(x,y))$ defines an inner product $\langle\cdot,\cdot\rangle_y^F$ on $T_xM$ which depends on the nonzero base vector $y$, and can be expressed as:
$$\langle u,v\rangle_y^F=g_{ij}(x,y)u^i v^j=\frac12\frac{\partial^2}{\partial s\partial t}F^2(x,y+su+tv)|_{s=t=0},$$
where
$u=u^i\partial_{x^i}$ and $v=v^j\partial_{x^j}$  are elements in $T_xM$.

The restriction of a Finsler metric $F$ to a tangent space is called a
{\it Minkowski norm}. Note that the notion of  a Minkowski norm $F$ satisfying the conditions (1)-(3) above can be defined on an arbitrary real vector space $\mathbf{V}$. The pair $(\mathbf{V},F)$ is called a {\it Minkowski space}.

For example, a Riemannian metric is a special Finsler metric $F$ such that
the Hessian matrix $(g_{ij}(x,y))$ is only relevant to $x$ on any standard
local coordinates. In this case, we usually refer the metric to be the globally
defined smooth section $g_{ij}(x)dx^idx^j$ of $\mathrm{Sym}^2(T^*M)$. A Randers metric is of the form $F=\alpha+\beta$, where $\alpha$ is a Riemannian metric and $\beta$ is a one-form satisfying some appropriate conditions (see \cite{BCS2000}). Randers metrics are the  simplest and  and most important non-Riemannian Finsler metrics. They can be generalized to $(\alpha,\beta)$-metrics which are of the form $F=\alpha\phi(\beta/\alpha)$,
where $\phi$ is a positive smooth function, and $\alpha$ and $\beta$ are similar as in the case of Randers metrics.

\subsection{Geodesic and geodesic spray}

A Finsler metric $F$ on a smooth manifold $M$ defines a possibly irreversible distance function $d_F(\cdot,\cdot)$ on  $M$. Then a geodesic can be defined as a piecewise smooth curve satisfying the locally minimizing principle. Note that on a Finsler space a geodesic need not be of constant speed. However, for the convenience  we will always parametrize a geodesic $c(t)$ to have positive constant speed, i.e.,  $F(\dot{c}(t))=F(\frac{d}{dt}c(t))=\mathrm{const}$. A geodesic in this sense  can be equivalently characterized as follows.

 Recall that {\it the geodesic spray} of a Finsler space $(M, F)$ is  a globally defined smooth vector field $\mathbf{G}$ on $TM\backslash 0$. On  a standard
local coordinates, it can be expressed as
\begin{equation}
\mathbf{G}=y^i\partial_{x^i}-2\mathbf{G}^i\partial_{y^i},
\end{equation}
where
\begin{equation}
\mathbf{G}^i=\frac{1}{4} g^{il}([F^2]_{x^k y^l}y^k - [F^2]_{x^l}).
\end{equation}
Then a curve $c(t)$ on $M$ is a geodesic of positive constant speed if and only if $(c(t),\dot{c}(t))$ is an integration curve
of $\mathbf{G}$. Thus on a standard local coordinates, a geodesic $c(t)=(c^i(t))$ satisfies the equations
\begin{equation}
\ddot{c}^i(t)+2 \mathbf{G}^i(c(t),\dot{c}(t))=0.
\end{equation}

\subsection{Flag curvature and S-curvature}

Curvature is the most important concept in Finsler geometry.
Some types of curvature are generalized from Riemannian geometry, which tell us how the space curves. We call them {\it Riemannian curvature}. Some others only appear for non-Riemannian Finsler spaces, that is,  they vanish on Riemannian manifolds. They tells us how different these spaces are from the Riemannian ones. We call them {\it non-Riemannian curvatures}.

Here we recall the definitions of  flag curvature and S-curvature, which are
two of the most important curvatures in Finsler geometry and are relevant to the main topics of this paper.

From the variational theory for a constant speed geodesic on a Finsler space $(M,F)$, one can deduce a similar Jacobi field equation as in Riemannian geometry, in which
there is a linear endomorphism $R_y^F:T_xM\rightarrow T_xM$ for any nonzero vector $y\in T_xM$.  This is   called the
{\it Riemann curvature}.
On a standard local coordinate system,
the Riemannian curvature can be presented as
$R^F_y=R^i_k(y)\partial_{x^i}\otimes dx^k:T_x M\rightarrow T_x M$, where
\begin{equation}
R^i_k(y)=2\frac{\partial}{\partial x^k}\mathbf{G}^i-y^j\frac{\partial^2}{\partial{x^j}\partial y^k}\mathbf{G}^i+2\mathbf{G}^j\frac{\partial^2}{\partial{y^j}\partial {y^k}}\mathbf{G}^i
-\frac{\partial}{\partial y^j}\mathbf{G}^i\frac{\partial}{\partial y^k}\mathbf{G}^j.
\end{equation}

Using the Riemannian curvature, the sectional curvature can be generalized to Finsler
geometry, which is called the {\it flag curvature}. Gievn $x\in M$, let $\mathbf{P}$ be a tangent plane in some $T_xM$, and $y$ a nonzero vector in $\mathbf{P}$.
Suppose $\mathbf{P}$ is linearly spanned by $y$ and $v$. Then the flag curvature
for the triple $(x,y,\mathbf{P})$ is defined as
\begin{equation}
K^F(x,y,\mathbf{P})=\frac{\langle R^F_y (v),v\rangle^F_y}
{\langle y,y\rangle^F_y \langle v,v\rangle^F_y-(\langle y,v\rangle^F_y)^2}.
\end{equation}
When $F$ is a Riemannian metric, it is just the sectional curvature for $(x,\mathbf{P})$, which is independent of the choice of $y$.

Z. Shen defines an important non-Riemannian curvature using the geodesic spray. It is  called {\it S-curvature} \cite{Sh97} in the literature.

Let $(M, F)$ be  a Finsler space and   $(x^1, x^2,\cdots, x^n, y^1,\cdots, y^n)$ a standard local coordinate system. The Busemann-Hausdorff volume form can be
defined as $dV_{BH}=\sigma(x)dx^1\cdots dx^n$, where
$$
\sigma(x)=\frac{\omega_n}{\mbox{Vol}\{(y^i)\in\mathbb{R}^n|F(x,y^i\partial_{x^i})<1\}},
$$
here $\mbox{Vol}$ denotes the volume of a subset with respect to the standard Euclidian metric on
$\mathbb{R}^n$, and $\omega_n=\mbox{Vol}(B_n(1))$. It is easily seen that the Busemann-Hausdorff
form is globally defined and does not depend on the specific coordinate system. On the other hand,
although the coefficient function $\sigma(x)$ is only locally defined and
depends on the choice of local coordinates $x=(x^i)$, the distortion function
\begin{equation*}
\tau(x,y)=\ln\frac{\sqrt{\det(g_{ij}(x,y))}}{\sigma(x)}
\end{equation*}
on $TM\backslash 0$ is independent of the local coordinates and globally defined.
The S-curvature $S(x,y)$ on $TM\backslash 0$ is  defined as the derivative of
$\tau(x,y)$ in the direction of the geodesic spray $\mathbf{G}(x,y)$.

\section{The positive curvature problem in Riemannian geometry}

In this section we collect the main results on the classification of positively curved spaces in   Riemannian geometry. It has been a very active topic  to find new examples of compact positively curved Riemannian manifolds. Generally speaking, this is
a very hard problem, because very few obstacles exist between those  with positive sectional curvature and with non-negative curvature. A remarkable open problem is the Hopf conjecture, which asks whether $S^2\times S^2$ admits a positively curved Riemannian metric.

On the other hand, imposing non-trivial Lie group actions as isometries on
the Riemannian manifolds can significantly reduce the complexity of the problem, which enables people to find new positively curved examples.
For example, positively curved Riemannian homogeneous spaces are totally
classified in a series of classical works during the 1960's and 1970's
(\cite{AW75, BB76, Ber61, Wallach1972}), which
can be summarized as the following classification theorem.
\begin{theorem}\label{overall-classification-theorem-in-Riemannian-case}
Up to local isometries, all positively curved Riemannian homogeneous spaces
belong to one of the following spaces:
\begin{description}
\item{\rm (1)} Compact rank one symmetric spaces
$S^{n-1}=\mathrm{SO}(n)/\mathrm{SO}(n-1)$,
$\mathbb{C}\mathrm{P}^{n-1}=\mathrm{SU}(n)/$ $\mathrm{S}(\mathrm{U}(n-1)\mathrm{U}(1))$,
$\mathbb{H}\mathrm{P}^{n-1}=\mathrm{Sp}(n)/\mathrm{Sp}(n-1)\mathrm{Sp}(1)$,
and $\mathrm{F}_4/\mathrm{Spin}(9)$.
\item{\rm (2)} Other homogeneous spheres and complex projective spaces, i.e.,
$S^{2n-1}=\mathrm{SU}(n)/$ $\mathrm{SU}(n-1)=\mathrm{U}(n)/\mathrm{U}(n-1)$,
$S^{4n-1}=\mathrm{Sp}(n)/\mathrm{Sp}(n-1)=\mathrm{Sp}(n)\mathrm{U}(1)/
\mathrm{Sp}(n-1)\mathrm{U}(1)=\mathrm{Sp}(n)\mathrm{Sp}(1)/\mathrm{Sp}(n-1)
\mathrm{Sp}(1)$,
$S^6=\mathrm{G}_2/\mathrm{SU}(3)$,
$S^7=\mathrm{Spin}(7)/\mathrm{G}_2$,
$S^{15}=\mathrm{Spin}(9)/\mathrm{Spin}(7)$, and
$\mathbb{C}\mathrm{P}^{2n-1}=\mathrm{Sp}(n)/\mathrm{Sp}(n-1)\mathrm{U}(1)$.
\item{\rm (3)} Berger's spaces $\mathrm{Sp}(2)/\mathrm{SU}(2)$ and
$\mathrm{SU}(5)/\mathrm{Sp}(2)S^1$.
\item{\rm (4)} Wallach's spaces $\mathrm{SU}(3)/T^2$, $\mathrm{Sp}(3)/\mathrm{Sp}(1)\mathrm{Sp}(1)\mathrm{Sp}(1)$ and
    $\mathrm{F}_4/\mathrm{Spin}(8)$, where $T^2$ and $\mathrm{Sp}(1)\mathrm{Sp}(1)\mathrm{Sp}(1)$ are the subgroups of all diagonal matrices.
\item{\rm (5)} Aloff-Wallach's spaces. These are $\mathrm{SU}(3)$-homogeneous
of the form  $$S_{k,l}=\mathrm{SU}(3)/\mathrm{diag}(z^k,z^l,z^{-k-l}), \, z\in\mathbb{C},$$ where the integers $k$ and $l$
    satisfy $kl(k+l)\neq 0$. They also admit $\mathrm{U}(3)$-homogeneous
    presentations. In particular, $S_{1,1}$
    can be expressed as $\mathrm{SU}(3)\times \mathrm{SO}(3)/\mathrm{U}(2)$.
\end{description}
\end{theorem}

Here are some historical remarks on the list in Theorem \ref{overall-classification-theorem-in-Riemannian-case}.
In 1961, M. Berger gave a  classification of positively curved Riemannian normal homogeneous spaces. His result is listed in  (1)-(3), in which the two spaces in (3) were  shown to admit positively curved normal Riemannian metrics, for the first time. However,  Berger missed $S_{1,1}=\mathrm{SU}(3)\times\mathrm{SO}(3)/\mathrm{U}(2)$ which also admits
positively curved Riemannian normal homogeneous metrics. This  was point out by B. Wilking in 1999 \cite{Wi1999}. Therefore in the following we will also refer
 $S_{1,1}=\mathrm{SU}(3)\times\mathrm{SO}(3)/\mathrm{U}(2)$ as
{\it the Wilking's space}.

In 1972, N. Wallach's achieved a classification  of even dimensional positively curved homogeneous spaces. His list consists of those even dimensional ones in (1) and (2), and the three new examples in (4) \cite{Wallach1972}.

In 1975,  S. Aloff and N. Wallach found the infinite sequence (5) of positive curved seven dimensional homogeneous spaces.  They conjectured that all positively
curved Riemannian homogeneous spaces had been found, i.e., are listed in Theorem \ref{overall-classification-theorem-in-Riemannian-case}. The conjecture was
proved by L. B\'{e}rard-Bergery in 1976 \cite{BB76}, in which he only needed to discuss the odd dimensional case.

In 2015, M. Xu and J. A. Wolf found a gap in \cite{BB76} in the argument for that $\mathrm{SO}(5)/\mathrm{SO}(2)=\mathrm{Sp}(2)/\mathrm{U}(1)$ can not be positively curved \cite{XW2015}. They proved that $\mathrm{SO}(5)/\mathrm{SO}(2)$ admits a Riemannian homogeneous metric which has positive sectional curvature for any linearly independent commuting pairs, so the algebraic method in \cite{BB76} does not work for this case. This phenomenon is very rare. The only other possible example is $\mathrm{Sp}(3)/\mathrm{diag}(z,z,q)$
with $z\in\mathbb{C}$ and $q\in\mathbb{H}$ \cite{WZ2015}\cite{XZi2016}. The gap in \cite{BB76} was fixed by a more analytic
argument by Wilking and Ziller. See \cite{WZ2015} and \cite{XW2015}.

The correctness of Theorem \ref{overall-classification-theorem-in-Riemannian-case} can be re-verified
by two approaches. One was provided by B. Wilking and W. Ziller \cite{WZ2015}, in which they used the classical fixed point set technique and some new methods from \cite{Wi2006}. Another was implied by the topic in this paper, i.e.,
in the process of  studying the classification problem of positively curved homogeneous Finsler spaces,  we can reduce the problem  to the Riemannian case (see \cite{XZi2016} or Theorem \ref{classification-reversible-alpha-beta}).

\section{Some techniques}
\label{section-3}
In this section  we present some techniques for the study of   positively curved homogeneous Finsler spaces.

\subsection{The Bonnet-Myers theorem and Synge's theorem}
The Bonnet-Myers theorem and Synge theorem's are valid in Finsler geometry \cite{BCS2000}.
By the Bonnet-Myers theorem, a complete connected
Finsler space with its Ricci scalar bounded below by a positive constant
must be compact. In particular, a homogeneous Finsler space with positive
flag curvature is compact.
On the other hand, by
Synge's theorem,  a compact connected positively curved Finsler space with odd dimension must be orientable, and the    fundamental group of a positively curved even-dimensional Finsler spaces is either trivial  or
$\mathbb{Z}_2$, according as it is orientable or not.

The proofs of the above  theorems involve the variational theory, i.e.,  the second derivative of the energy functional at a geodesic. For the Bonnet-Myers theorem, the end points of the geodesic are unmoved during the variation. For  Synge's theorem, each curve is closed during the variation. In both cases, the variation theory is very similar to the Riemannian context \cite{BCS2000}. Notice that there is
a third variational theory for the energy functional at a geodesic, which moves the end points along two geodesics (or more generally two totally geodesic
submanifolds, as for Frankel's theorem). It can not be easily generalized to
Finsler geometry because of the change of base vectors. This fact interprets why  Frankel's theorem and some  related works have not been generalized from
Riemannian geometry to Finsler geometry.

By the Bonnet-Myers Theorem,   a homogeneous Finsler space $(M,F)$ with positive flag curvature can be written as a coset space $M=G/H$, where  $G$ is a compact Lie group. In some cases,
we need to assume  $G$ to be simply connected. If in this case $G$   has a center of positive dimension, then we can assume $G$ to be quasi-compact, i.e.,  $\mathfrak{g}=\mathrm{Lie}G$ is
a compact Lie algebra which admits an $\mathrm{Ad}(G)$-invariant inner product. This is  the start point  for the linear and algebraic setups for
classifying positively curved homogeneous Finsler spaces.

\subsection{Totally geodesic subspace and fixed point set technique}
A {\it totally geodesic subspace} $(N,F|_N)$ of a Finsler space $(M,F)$ is an immersed submanifold $N$ endowed with the subspace metric $F|_N$, such that each geodesic of $(N,F|_N)$ is also a geodesic of $(M,F)$. It is not  hard to verify the equalities $$R^{F|_N}_y=R^F_y\mbox{ and }
K^{F|_N}(x,y,\mathbf{P})=K^F(x,y,\mathbf{P})$$
for any $x\in N$, tangent plane $\mathbf{P}$ in $T_xN$ and nonzero tangent vector $y\in \mathbf{P}$.

The fixed point set $\mathrm{Fix}(\mathcal{A})$ of a family $\mathcal{A}$ of isometries of a Finsler space $(M,F)$
is totally geodesic \cite{De2008}. If  $(M,F)$ is homogeneous, then it is also true for   any connected component $\mathrm{Fix}(\mathcal{A})_o$ of $\mathrm{Fix}(\mathcal{A})$.
 This observation leads to an important fact that if  $(M,F)$ is a connected positively curved homogeneous Finsler space and $\dim \mathrm{Fix}(\mathcal{A})_o>1$,
then $\mathrm{Fix}(\mathcal{A})_o$ is also a connected positively curved homogeneous Finsler space. In the following, the above  method or technique
will be indicated as {\it totally geodesic technique} or
 {\it fixed point set technique}.

One application of this technique is to consider $\mathrm{Fix}(T_H)_o$ in a positively curved homogeneous Finsler space
$M=G/H$, where $G$ is compact or quasi-compact, $H$ is the isotropy subgroup at $o=eH$, and $T_H$ is a maximal torus in $H$. Then the component
$(\mathrm{Fix}(T_H)_o, F|_{\mathrm{Fix}(T_H)_o})$ is locally isometric to $(G',F')$ where $G'$ is a compact Lie group and $F'$ is a left invariant
Finsler metric. So we have the equivalence between the two statements in
the following theorem.

\begin{theorem} \label{rank-inequ-thm}
(1) If a compact Lie group $G$ admits a left invariant
positively curved Finsler metric, then it must be isomorphic either to $\mathrm{SU}(2)$
or to $\mathrm{SO}(3)$.

(2) Suppose a compact coset space $G/H$,  where $G$ is a quasi compact Lie group,
admits a positively curved $G$-homogeneous Finsler metric. Then
we have $\mathrm{rank}\mathfrak{g}\leq\mathrm{rank}\mathfrak{h}+1$. More precisely, if
$\dim G/H$ is even, then $\mathrm{rank}\mathfrak{g}=\mathrm{rank}\mathfrak{h}$;
 if $\dim G/H$ is odd, then $\mathrm{rank}\mathfrak{g}=\mathrm{rank}\mathfrak{h}+1$.
\end{theorem}

 By Theorem \ref{rank-inequ-thm}, to prove Theorem \ref{overall-classification-theorem-in-Riemannian-case}, we only need to prove the first statement.
There are three methods meeting our goal. The first was given by S. Deng and Z. Hu \cite{DH2013},
using a similar argument as N. Wallach in \cite{Wallach1972}. The second was suggested by W. Ziller, using the totally geodesic technique ( see the comment on the rank inequality in Section 2 of \cite{XZi2016}). The third was a purely computational proof found by L. Huang, using the technique of B. Wilking (see Lemma 1.1 in \cite{WZ2015} or Lemma 5.2 in \cite{XW2015}) and his homogeneous flag curvature formula \cite{Huang2013}.

Another application of the totally geodesic technique is to consider
$\mathrm{Fix}(\iota)_o$ of
a noncentral involution $\iota\in H$, i.e.,  $G^\iota=C_G(\iota)\neq G$ and $\iota^2=e\in G$. All possible cases of $G^\iota$ and
$\mathrm{Fix}(\iota)_o=G^\iota/H^\iota$ are
given by the classification theory for symmetric spaces and the induction method. The classification can be finally reduced to a very short list of
the cases treated in \cite{WZ2015} and \cite{XZi2016}.

\subsection{Finsler submersion and homogeneous flag curvature formula}
Submersions can be defined between Finsler spaces \cite{PD2001}. A smooth map between two Finsler spaces $\pi:(M_1,F_1)\rightarrow (M_2,F_2)$ is called a {\it Finsler
submersion} if its tangent map at each point $\pi_*:(T_xM_1,F_1(x,\cdot))
\rightarrow(T_{\pi(x)}M_2,F_2(\pi(x),\cdot))$ is a {\it linear Finsler submersion}, i.e. $\pi_*$ maps the $F_1$-unit disk in $T_xM_1$ onto the $F_2$-unit disk in $T_{\pi(x)}M_2$.

Horizonal lifting is a crucial concept
for Finsler submersions. For example, given a Finsler submersion $\pi:(M_1,F_1)\rightarrow (M_2,F_2)$, and a
 triple $(\pi(x),y,\mathbf{P})$, where $\mathbf{P}$ is a tangent plane in
$T_{\pi(x)}M$ and $y$ is a nonzero vector in $\mathbf{P}$, we have a  horizonal
lifting $(x,\tilde{y},\tilde{\mathbf{P}})$ of $(\pi(x),y,\mathbf{P})$. The  important point is that  we have an inequality between the  flag curvatures of them \cite{PD2001}
\begin{equation}\label{submersion-flag-curvature-inequality}
K^{F_1}(x,\tilde{y},\tilde{\mathbf{P}})\leq K^{F_2}(\pi(x),y,\mathbf{P}).
\end{equation}

Finsler submersion can be applied to the study of homogeneous Finsler geometry, in that for any homogeneous Finsler space $(G/H,F)$, we can find
a Finsler metric $\bar{F}$ on $G$ such that $\bar{F}$ is left $G$-invariant and right $H$-invariant, and the canonical projection $\pi:(G,\bar{F})\rightarrow(G/H,F)$ is a Finsler submersion. Conversely,
given  a Finsler metric $\bar{F}$ on $G$,  a homogeneous Finsler metric $F$
can be induced by $\bar{F}$ using the submersion method.
We generally call this kind of   methods   the {\it Finsler submersion technique}.

This technique has been applied in our classification of positively curved normal homogeneous
Finsler spaces. A homogeneous Finsler space $(G/H,F)$ is call {\it normal} or {\it $G$-normal}, if it is induced by a bi-invariant Finsler metric $\bar{F}$ on $G$
such that the canonical projection $\pi:(G,\bar{F})\rightarrow (G/H,F)$ is
a Finsler submersion.

A significant achievement of the Finsler submersion technique is  the following theorem \cite{XDHH2014}.
\begin{theorem}\label{homogeneous-flag-curvature-formula-thm}
Let $(G/H,F)$ be a connected homogeneous Finsler space, and $\mathfrak{g}=\mathfrak{h}+\mathfrak{m}$ be an $\mathrm{Ad}(H)$-invariant
decomposition for $G/H$. Then for any linearly independent commuting pair
$u$ and $v$ in $\mathfrak{m}$ satisfying $\langle [u,\mathfrak{m}],u\rangle_u^F=0$, we have
\begin{equation}\label{homogeneous-flag-curvature-formula}
K^F(o,u,u\wedge v)=\frac{\langle U(u,v),U(u,v)\rangle_u^F}{\langle u,u\rangle_u^F\langle v,v\rangle_u^F-[\langle u,v\rangle_u^F]^2},
\end{equation}
where $U(u,v)$ is the map from $\mathfrak{m}\times\mathfrak{m}$
to $\mathfrak{m}$ determined by
$$\langle U(u,v),w\rangle_u^F=\frac12(\langle[w,u]_\mathfrak{m},v\rangle_u^F+
\langle[w,v],u\rangle_u^F),\quad\forall w\in\mathfrak{m},$$
where $[\cdot,\cdot]_\mathfrak{m}=\mathrm{pr}_m\circ[\cdot,\cdot]$ and
$\mathrm{pr}_\mathfrak{m}$ is the projection with respect to the given
$\mathrm{Ad}(H)$-invariant decomposition.
\end{theorem}

The equality (\ref{homogeneous-flag-curvature-formula}) is called a
{\it homogeneous flag curvature formula}.
Though it is not presented in the most general form (it can
be deduced from the more complicated flag curvature formula of L. Huang \cite{Huang2013}),  its simpleness and practicability mark it a radiant star in
homogeneous Finsler geometry.

The homogeneous flag curvature formula (\ref{homogeneous-flag-curvature-formula}) is crucial  when we need to
prove that some special compact coset spaces $G/H$ does not admit positive flag curvature. It is usually applied as the following. Assume conversely that $G/H$ does admit a positively curved Finsler metric $F$. Then we just need to find a linearly
independent commuting pair $u$ and $v$ in $\mathfrak{m}$ such that
$$\langle[u,\mathfrak{m}]_\mathfrak{m},u\rangle_u^F=
\langle[v,\mathfrak{m}]_\mathfrak{m},u\rangle_u^F=
\langle[u,\mathfrak{m}]_\mathfrak{m},v\rangle_u^F=0.$$
Then Theorem \ref{homogeneous-flag-curvature-formula-thm} implies
that $K^F(o,u,u\wedge v)=0$, which is a contradiction.

For most cases, the argument that $G/H$ does not admit positive flag curvature
has a similar pattern. This method has been summarized  as several key lemmas (see \cite{XD2016} or, Lemma \ref{key-lemma} and Lemma \ref{key-lemma-2}). Using those key lemmas, the classification of positively curved homogeneous Finsler
spaces can be reduced to an algebraic problem.

\subsection{Linear and algebraic setups for a positively curved homogeneous Finsler space}

 Let $G/H$ be a positively curved homogeneous Finsler space, where $G$ is a quasi-compact group. Fix an $\mathrm{Ad}(G)$-invariant inner product
$\langle\cdot,\cdot\rangle_{\mathrm{bi}}$ on $\mathfrak{g}=\mathrm{Lie}(G)$
and an orthogonal decomposition
$\mathfrak{g}=\mathfrak{h}+\mathfrak{m}$, where
$\mathfrak{h}=\mathrm{Lie}(H)$. Denote $\mathrm{pr}_\mathfrak{h}$
and $\mathfrak{pr}_\mathfrak{m}$ the orthogonal projections to $\mathfrak{h}$
and $\mathfrak{m}$, respectively.

Fix a Cartan subalgebra $\mathfrak{t}$ of $\mathfrak{g}$ such that
$\mathfrak{t}\cap\mathfrak{h}$ is a Cartan subalgebra of $\mathfrak{h}$.
With respect to $\mathfrak{t}$ and $\mathfrak{t}\cap\mathfrak{h}$,
we have the following decompositions of $\mathfrak{g}$ and $\mathfrak{h}$ respectively:
\begin{eqnarray}
\mathfrak{g}&=&\mathfrak{t}+\sum_{\alpha\in\Delta(\mathfrak{g},\mathfrak{t})}
\mathfrak{g}_{\pm\alpha},\label{decomp-g}\\
\mathfrak{h}&=&\mathfrak{t}\cap\mathfrak{h}+
\sum_{\alpha'\in\Delta(\mathfrak{h},\mathfrak{t}\cap\mathfrak{h})}
\mathfrak{h}_{\pm\alpha'},\label{decomp-h}
\end{eqnarray}
where $\Delta(\mathfrak{g},\mathfrak{t})$ and
$\Delta(\mathfrak{h},\mathfrak{t}\cap\mathfrak{h})$ are the root systems, and $\mathfrak{g}_{\pm\alpha}$ and
$\mathfrak{h}_{\pm\alpha'}$ are the root planes of $\mathfrak{g}$ and $\mathfrak{h}$, respectively. Notice that in our conventions, $\Delta(\mathfrak{g},\mathfrak{t})$ (resp.
$\Delta(\mathfrak{h},\mathfrak{t}\cap\mathfrak{h})$) is viewed as a subsets of
$\mathfrak{t}$ (resp. $\mathfrak{t}\cap\mathfrak{h}$) through the inner product $\langle\cdot,\cdot\rangle_{\mathrm{bi}}$.

 If $\dim G/H$ is even, then the rank inequality implies that $\mathfrak{t}\subset\mathfrak{h}$. Then a  root plane of $\mathfrak{g}$ is contained either  in $\mathfrak{h}$ or $\mathfrak{m}$. Thus $H$ a {\it regular} subgroup of $G$, i.e.,  each root plane of $\mathfrak{h}$ is also a root
plane of $\mathfrak{g}$.

If  $\dim G/H$ is odd, then the rank inequality implies that
$\dim\mathfrak{t}\cap\mathfrak{m}=1$. In this case, the decompositions (\ref{decomp-g}) and
(\ref{decomp-h}) are compatible with $\mathfrak{g}=\mathfrak{h}+\mathfrak{m}$,
 namely, given a  nonzero $\alpha'\in\mathfrak{t}\cap\mathfrak{h}$, denoting $\hat{\mathfrak{g}}_{\pm\alpha'}=\sum_{\mathrm{pr}_\mathfrak{h}
(\alpha)=\alpha'}\mathfrak{g}_{\pm\alpha}$,  we have
\begin{eqnarray*}
\hat{\mathfrak{g}}_{\pm\alpha'}=\hat{\mathfrak{g}}_{\pm\alpha'}\cap\mathfrak{h}+
\hat{\mathfrak{g}}_{\pm\alpha'}\cap\mathfrak{m}.
\end{eqnarray*}

If  $\alpha'$ is a root of $\mathfrak{h}$,
then we have $\hat{\mathfrak{g}}_{\pm\alpha'}\cap\mathfrak{h}=\mathfrak{h}_{\pm\alpha'}$.
Otherwise $\hat{\mathfrak{g}}_{\pm\alpha'}\subset\mathfrak{m}$. So
each root of $\mathfrak{h}$ must be equal to $\mathrm{pr}_\mathfrak{h}(\alpha)$
for some root of $\mathfrak{g}$. If a root plane $\mathfrak{h}_{\pm\alpha'}$
of $\mathfrak{h}$ is not a root plane of $\mathfrak{g}$, there must be
a linear independent pair of roots $\alpha$ and $\beta$  of $\mathfrak{g}$,
such that $\mathrm{pr}_\mathfrak{h}(\alpha)=\mathrm{pr}_\mathfrak{h}(\beta)
=\alpha'$.

To summarize, we can divide the odd dimensional positively curved homogeneous Finsler spaces $(G/H,F)$ into three categories:
\begin{description}
\item{\bf Case I.} $H$ is a regular subgroup of $G$, i.e.,  each root plane of $\mathfrak{h}$ is a root plane of $\mathfrak{g}$.
\item{\bf Case II.} $H$ is not a regular subgroup of $G$, and there exists a root $\alpha'$ of $\mathfrak{h}$, and a linearly independent pair of roots $\alpha$ and $\beta$ from different simple factors of $\mathfrak{g}$, such that
    $\mathrm{pr}_\mathfrak{h}(\alpha)=\mathrm{pr}_\mathfrak{h}(\beta)
=\alpha'$.
\item{\bf Case III.} $H$ is not a regular subgroup of $G$, and there exists a root $\alpha'$ of $\mathfrak{h}$, and a linearly independent pair of roots $\alpha$ and $\beta$ from the same simple factor of $\mathfrak{g}$, such that
    $\mathrm{pr}_\mathfrak{h}(\alpha)=\mathrm{pr}_\mathfrak{h}(\beta)
=\alpha'$.
\end{description}
 Therefore, to classify positively curved homogeneous Finsler spaces, we just need to find the classification list for Case I-III.

\section{Classification of positively curved homogeneous Finsler spaces}

\subsection{Introduction to the main problem}
The following important problem is of great interest to Finsler geometers.
\begin{problem}\label{main-prob}
Classify all compact coset spaces which admit homogeneous Finsler metrics with
positive flag curvature.
\end{problem}
This is the key step in the  study of  Finsler spaces with positive flag curvature.
More importantly, a complete classification will help us understand the flag curvature in Finsler geometry.

Though flag curvature shares many important properties with sectional curvature as stated in the previous section, its formulas by local coordinates are much more complicated than in Riemannian
geometry. Until recent years, Problem \ref{main-prob} in the general sense has not been touched, except for a few very special cases \cite{DH2013} \cite{HD11}.

Building the linear and algebraic
setups for a positively curved homogeneous Finsler space $(G/H,F)$ marks the
first milestone for studying this classification problem in the general sense.
It was established when
we classified  normal homogeneous Finsler spaces in \cite{XD-Normal-homogeneous-Finsler-spaces}, where we  got the same classification list as M. Berger's in \cite{Ber61}. Notice that we also established some other important techniques in \cite{XD-Normal-homogeneous-Finsler-spaces}, such as  the totally geodesic technique and the Finsler submersion technique.

Finding the homogeneous flag curvature formula (\ref{homogeneous-flag-curvature-formula})
marks an even more splendid progress for studying Problem \ref{main-prob} \cite{XDHH2014}. The formula itself is surprisingly beautiful, powerful and simple.
Due to this discovery, the algebraic method for
classifying positively curved Riemannian homogeneous spaces can be reformulated in the Finsler context. In the even dimensional case, positively curved homogeneous Finsler spaces can be completely classified, from which we get the
same classification list as N. Wallach in the Riemannian context \cite{XDHH2014}. In the odd dimensional case, we  need to assume the metric to be reversible. Then
the homogeneous flag curvature formula can be effectively applied. Based on this idea,  we generalized the classification work of L. B\'{e}rard-Bergery to Finsler geometry, and showed that an odd dimensional positively curved reversible
homogeneous Finsler space either admits a positively curved Riemannian homogeneous metric, or  belongs to a short list of five exceptional cases in \cite{XZi2016}.

We will present more details for these classification works in
Subsection \ref{g-c-t}. Since in  these works  Problem \ref{main-prob} is actually considered
in a general sense, we will refer  to the related results as {\it general classifications}.

Notice that when the dimension is odd and the metric is reversible, the classification for the hardest case, i.e.,  Case I, has not been completely finished.
In fact, all the five undetermined cases belong to Case I. When the dimension is odd and
the metric is irreversible,  even less is known for the general classifications.

 This fact led us to the study of  the \textit{unsolved parts} of Problem \ref{main-prob}, imposing some special recipe in Finsler geometry. In particular, in some of our previous works we either assume that  the metrics belong to some special class, such as  Randers
metrics,   $(\alpha,\beta)$-metrics, or add more curvature conditions, e.g.,  spaces with positive falg curvature and vanishing S-curvature. This consideration will be referred to as
{\it special classifications}.

In establishing  the special classifications, some special methods can be applied besides those listed in Section 3. For non-Riemannian homogeneous
Randers spaces and $(\alpha,\beta)$-spaces with positive flag curvature and
vanishing S-curvature, we can get a complete classification. Notice that these
homogeneous Finsler spaces are generally irreversible and they belong to Case I. For non-Riemannian reversible homogeneous $(\alpha,\beta)$-spaces with positive flag curvature, the classification list is still incomplete, but
the number of undetermined cases can be reduced to $2$. We will give more details for these special classifications in Subsection \ref{s-c-t}.

To conclude this subsection, we give two remarks.

First, all the classifications
mentioned above only generalize the algebraic methods from the Riemannian context which exclude unwanted spaces from the list of possible candidates.
So they are not sufficient to give a complete classification list because of
the two undetermined cases
\begin{eqnarray*}
\mathrm{SO}(5)/\mathrm{SO}(2)&=&\mathrm{Sp}(2)/\mathrm{diag}(z,z)\mbox{ with }z\in\mathbb{C},\\
\mathrm{Sp}(3)/\mathrm{Sp}(1)_{(3)}S^1_{(1,1,0)}&=&
\mathrm{Sp}(3)/\mathrm{diag}(z,z,q)\mbox{ with }z\in\mathbb{C}\mbox{ and }q\in\mathbb{H}.
\end{eqnarray*}
The reason is explained in
\cite{WZ2015} \cite{XZh2016},  and sketchily explained in Section 2.

On the other hand, due to the complexity of calculation in Finsler geometry,
we have no clue on how to generalize the analytic methods in the Riemannian context, so we can not find
new positively curved examples by these classifications.

Besides the classifications in the space level, the classifications in the metric level, i.e.,  giving complete and explicit descriptions on the space of all positively curved homogeneous Finsler metrics, is another important research project. For example, in Riemannian geometry, it has been thoroughly discussed for homogeneous spheres \cite{VZ2009}. In Finsler geometry, this is a much harder problem than the classification of the spaces. In most occasions, we can only
 determine whether  a positively curved Riemannian homogeneous metric admits generic non-Riemannian perturbations, which will lead to positively curved non-Riemannian homogeneous Finsler metrics. Up to now, the
only complete classification achieved in the metric level is for  non-Riemannian homogeneous Randers metrics with positive flag curvature and vanishing S-curvature \cite{HD11}.

\subsection{General classifications}
\label{g-c-t}

First we consider the classification of positively curved normal homogeneous
Finsler spaces. The definition of normal homogeneous Finsler spaces was first made clear in \cite{XD2015}. Recall that if $G$ is a quasi-compact Lie group and $G/H$ is a coset space of $G$, where $H$ is a closed subgroup of $G$, then a Finsler metric $F$ on $G$ is called normal if it is subdued by a bi-invariant Finsler metric $\bar{F}$ on $G$, which means that the differential map of the natural projection from $G$ onto $G/H$ maps the indicatrix of $(T_e(G), \bar{F})$ onto the indicatrix of $(T_o(G/H), F)$, where $o$ is the origin of $G/H$.

Let $(G/H,F)$ be a $G$-normal homogeneous Finsler space, where $G$ is a
quasi-compact Lie group.  Fix an $\mathrm{Ad}(G)$-invariant inner product
on $\mathfrak{g}$ and  the corresponding orthogonal decomposition
$\mathfrak{g}=\mathfrak{h}+\mathfrak{m}$.
A subalgebra $\mathfrak{s}$ of $\mathfrak{g}$ is called a {\it flat splitting subalgebra} (or {\it FSS} in short) with respect to the above settings if the following conditions are
satisfied:
\begin{description}
\item{\rm (1)} $\mathfrak{s}$ is the intersection of a family of
Cartan subalgebras of $\mathfrak{g}$.
\item{\rm (2)} $\mathfrak{s}=\mathfrak{s}\cap\mathfrak{h}+
\mathfrak{s}\cap\mathfrak{m}$.
\item{\rm (3)} $\dim\mathfrak{s}\cap\mathfrak{m}\geq 2$.
\end{description}
A special case of FSS is {\it flat splitting Cartan subalgebra} or {\it FSCS} in short, when it is a Cartan subalgebra itself.

If $\mathfrak{s}$ is a FSS for a compact
normal homogeneous Finsler $(G/H,F)$, then
$S=\exp\mathfrak{s}\cap\mathfrak{m}$ is totally geodesic  in $(G/H,F)$. Notice that it is a torus, and the induced  metric $F|_S$ is
left invariant, i.e.,  a flat metric on $S$. So the existence of a FSS provides
an obstacle for $(G/H,F)$ to have positive flag curvature.

Notice that for any
closed Lie group $K$   satisfying
$\mathfrak{h}\subset\mathfrak{k}\subset\mathfrak{g}$, where  $\mathfrak{k}=\mathrm{Lie} (K)$, the $G$-normal homogeneous
Finsler metric $F$ on $G/H$  induces a $G$-normal homogeneous
Finsler metric on $G/H_0$ which is locally isometric to $F$  (denoted by the same $F$),  and induces
a $G$-normal homogeneous metric $F'$ on $G/K$ such that the canonical projection map
$\pi:(G/H_0,F)\rightarrow (G/K,F')$ is a Finsler submersion. If $F$ is
positively curved and $\dim G/K>1$, then  $F'$ is also positively curved on $G/K$.

In summarizing, we get an  obstacle for homogeneous positive curvature, which can be   stated as in  the following theorem; See \cite{XD-Normal-homogeneous-Finsler-spaces}.

\begin{theorem} \label{FSS-obstacle-thm}
Let $F$ be a  normal homogeneous Finsler metric  on
$G/H$ with  positive flag curvature, induced by an $\mathrm{Ad}(G)$-invariant Minkowski norm on $\mathfrak{g}$. Let $K$ be a closed subgroup  of $G$ with $\mathfrak{h}\subset\mathfrak{k}\neq\mathfrak{g}$. Fix  an  orthogonal decomposition
$\mathfrak{g}=\mathfrak{k}+\mathfrak{p}$ with respect to $\langle\cdot,\cdot\rangle_{\mathrm{bi}}$ on $\mathfrak{g}$.
Then $G/K$ admits a positively
curved normal homogeneous Finsler metric, and there does not exist
any FSS for $G/K$ with respect to the orthogonal decomposition of an $\mathrm{Ad}(G)$-invariant inner product on $\mathfrak{g}$.
\end{theorem}

Theorem \ref{FSS-obstacle-thm} reduces the classification for positively
curved normal homogeneous Finsler spaces to a totally algebraic problem.
Based on this theorem and  a case-by-case discussion, we get  a complete classification result \cite{XD-Normal-homogeneous-Finsler-spaces}:

\begin{theorem} Let $G$ be a connected compact Lie group and $H$ a closed subgroup of $G$. Then there exists a $G$-invariant normal homogeneous Finsler
metric on $G/H$ with positive flag curvature if and only if there exists a
normal homogeneous Riemannian metric on $G/H$ with positive sectional curvature.
\end{theorem}

Therefore,  the complete classification list of positively curved normal homogeneous
Finsler spaces consists of (1)-(3) in Theorem \ref{overall-classification-theorem-in-Riemannian-case}. Notice that in the classification list given by Theorem 1.1 of \cite{XD-Normal-homogeneous-Finsler-spaces}, we missed the homogeneous complex
projective space $\mathbb{C}\mathrm{P}^{2n-1}=\mathrm{Sp}(n)/\mathrm{Sp}(n-1)\mathrm{U}(1)$, which was provided by the first paragraph in page 19 of \cite{XD-Normal-homogeneous-Finsler-spaces}.

This  method can also be applied to $\delta$-homogeneous Finsler
spaces \cite{XZh2016}. The
$\delta$-homogeneity has several equivalence definitions, one of which has
the same pattern as normal homogeneity. Any $G$-$\delta$-homogeneous Finsler
metric can be approximated by a sequence of $G$-normal homogeneous Finsler metrics. An FSS for a $G$-$\delta$-homogeneous Finsler space $(G/H,F)$ also
provides a totally geodesic flat subspace.  Then a similar statement as Theorem
\ref{FSS-obstacle-thm} and a case-by-case discussion show that
a compact coset space $G/H$ admits positively curved
$G$-$\delta$-homogeneous Finsler metrics if and only if it admits positively curved
$G$-normal homogeneous Finsler metrics. So the three classification lists are the same.

Next we turn to the classification of even dimensional positively curved homogeneous Finsler spaces.

In \cite{XDHH2014}, together with L. Huang and Z. Hu, we used the homogeneous flag curvature formula (\ref{homogeneous-flag-curvature-formula}) in Theorem \ref{homogeneous-flag-curvature-formula-thm} to prove the following lemma.

\begin{lemma} \label{strong-ortho-lemma}
Let $(G/H,F)$ be an even dimensional homogeneous Finsler space,  where $G$ is
a compact Lie group. Fix
a Cartan subalgebra $\mathfrak{t}$ in $\mathfrak{h}=\mathrm{Lie}(H)$. Then for
any linearly independent pair of roots $\alpha$ and $\beta$ of $\mathfrak{g}$ which are not roots of $\mathfrak{h}$, either $\alpha+\beta$ or $\alpha-\beta$
is  a root of $\mathfrak{g}$.
\end{lemma}

The compact coset spaces satisfying the property for roots in Lemma \ref{strong-ortho-lemma} was completely classified by N. Wallach \cite{Wallach1972},
which gives his classification  of even dimensional positively curved
Riemannian homogeneous spaces. So we got the same classification list
in the Finsler context,
i.e., we have  the following main theorem in \cite{XDHH2014}.

\begin{theorem} Let $G$ be a compact connected simply connected Lie group and $H$ a connected closed subgroup such that the dimension of the coset space
$G/H$ is even. Suppose that there exists a $G$-invariant Finsler metric on
$G/H$ with positive flag curvature. Then there exists a $G$-invariant Riemannian metric on $G/H$ with positive sectional curvature.
\end{theorem}

The complete classification list of even dimensional positively curved homogeneous Finsler spaces consists of the even dimensional ones in (1) and (2), and the three Wallach's spaces in (4) of the list in Theorem \ref{overall-classification-theorem-in-Riemannian-case}.

Finally, we consider the classification of odd dimensional positively curved
reversible homogeneous Finsler spaces.

The reason to assume the reversibility for the homogeneous Finsler metric $F$ is the following. When we apply the homogeneous flag curvature formula
(\ref{homogeneous-flag-curvature-formula}) to prove that $(G/H,F)$ is not positively curved, we need to find a linearly independent commutative pair $u$ and $v$
from $\mathfrak{m}$ such that
$$\langle[u,\mathfrak{m}]_\mathfrak{m},u\rangle_u^F=
\langle[v,\mathfrak{m}]_\mathfrak{m},u\rangle_u^F=
\langle[u,\mathfrak{m}]_\mathfrak{m},v\rangle_u^F=0.$$
 Since $\dim G/H$ is odd, $\mathfrak{t}\cap\mathfrak{m}$ is a one-dimensional subspace. To get the above equalities, we will need the condition
$\langle \mathfrak{t}\cap\mathfrak{m},u\rangle_u^F=0$. This is not true when $F$
is reversible, but incorrect when $F$ is not.

We now recall the following
two key lemmas in \cite{XD2016}.

\begin{lemma}\label{key-lemma}
Let $F$ be a positively curved homogeneous Finsler metric on the odd dimensional coset space $G/H$. If $\alpha$ is a root of $\mathfrak{g}$
contained in $\mathfrak{t}\cap\mathfrak{h}$, and it is the only root of
$\mathfrak{g}$ contained in $\alpha+\mathfrak{t}\cap\mathfrak{m}$, then it
must be a root of $\mathfrak{h}$ and we have
$\mathfrak{h}_{\pm\alpha}=\hat{\mathfrak{g}}_{\pm\alpha}
=\mathfrak{g}_{\pm\alpha}$.
\end{lemma}

\begin{lemma}\label{key-lemma-2} Let $F$ be a reversible positively curved homogeneous Finsler metric on an odd dimensional coset space $G/H$. Then there does not exist a pair of linearly independent roots $\alpha$ and $\beta$ of $\mathfrak{g}$ such that the followings hold simultaneously:
\begin{description}
\item{\rm (1)} Neither $\alpha$ nor $\beta$ is a root of $\mathfrak{h}$.
\item{\rm (2)} $\pm\alpha$ are the only roots of $\mathfrak{g}$ in
$\mathbb{R}\alpha+\mathfrak{t}\cap\mathfrak{m}$.
\item{\rm (3)} $\pm\beta$ are the only roots of $\mathfrak{g}$
in $\mathbb{R}\alpha\pm\beta+\mathfrak{t}\cap\mathfrak{m}$.
\end{description}
\end{lemma}

Both the key lemmas are proved by deducing a contradiction. We just need to find a suitable
linearly independent commutative pair $u$ and $v$ from $\mathfrak{m}$ such that $K^F(o,u,u\wedge v)=0$.  Notice that in Lemma \ref{key-lemma}, we do not need the reversibility, and in Lemma \ref{key-lemma-2}, the case that $\alpha\in\mathfrak{t}\cap\mathfrak{m}$ is
in fact covered by Lemma \ref{key-lemma}.

These two key lemmas, especially the second one, are crucial for classifying
odd dimensional positively curved reversible homogeneous Finsler spaces.
Using them, we can exclude many  coset spaces in the list of
possible candidates. Using this method, we classified all the coset spaces one by one in each subcase of Case II and Case III, which can be summarized as the following main theorem in \cite{XD2016}.

\begin{theorem}
Let $(G/H,F)$ be an odd-dimensional positively curved reversible homogeneous Finsler space. If $H$ is not regular in $G$, then
$G/H$ admits a $G$-invariant Riemannian metric with positive curvature.
\end{theorem}

To be precise, if $G/H$ is of Case II, then it is equivalent to one of
the following coset spaces:
\begin{description}
\item{\rm (1)} Homogeneous spheres $S^3=\mathrm{SO}(4)/\mathrm{SO}(3)$ and
$S^{4n-1}=\mathrm{Sp}(n)\mathrm{Sp}(1)/\mathrm{Sp}(n-1)\mathrm{Sp}(1)$ with $n>1$;
\item{\rm (2)} Wilking's space $\mathrm{SU}(3)\times\mathrm{SO}(3)/\mathrm{U}(2)$.
\end{description}
If  $G/H$ is of Case III, then it is equivalent to one of
the following coset spaces:
\begin{description}
\item{\rm (1)} Homogeneous spheres
$S^{2n-1}=\mathrm{SO}(2n)/\mathrm{SO}(2n-1)$ with $n>2$,
$S^7=\mathrm{Spin}(7)/\mathrm{G}_2$,
and
$S^{15}=\mathrm{Spin}(9)/\mathrm{Spin}(7)$;
\item{\rm (2)}
Berger's spaces $\mathrm{SU}(5)/\mathrm{Sp}(2)\mathrm{U}(1)$
and $\mathrm{Sp}(2)/\mathrm{SU}(2)$.
\end{description}

Among all the three cases, Case I is the hardest one. Up to now, the
general classification for this case is still incomplete.
In \cite{XD2016}, we used
the homogeneous flag curvature formula (\ref{homogeneous-flag-curvature-formula}) to discuss this case and proved another main theorem.

\begin{theorem}\label{Case-I-preparation-thm}
Let $(G/H,F)$ be an odd-dimensional positively curved reversible homogeneous Finsler space. If $H$ is a regular subgroup of $G$,  then
there are only the following two cases:
\begin{description}
\item{\rm (1)} $G/H$ is equivalent to the homogeneous spheres
$S^{2n-1}=\mathrm{U}(n)/\mathrm{U}(n-1)$,
$S^{4n-1}=\mathrm{Sp}(n)\mathrm{U}(1)/\mathrm{Sp}(n-1)\mathrm{U}(1)$ with $n>1$, or the $\mathrm{U}(3)$-homogeneous Aloff-Wallach's spaces.
\item{\rm (2)} $G/H$ is equivalent to an odd-dimensional reversible positively curved homogeneous Finsler space $G'/H'$ such that $G'$
    is a compact simple Lie group and $H'$ is a regular subgroup of $G'$.
\end{description}
\end{theorem}

Notice that the only relevant cases in the list of Theorem \ref{overall-classification-theorem-in-Riemannian-case} which do not appear in the above theorem are $S^{2n-1}=\mathrm{SU}(n)/\mathrm{SU}(n-1)$,
$S^{4n-1}=\mathrm{Sp}(n)/\mathrm{Sp}(n-1)$, and the
$\mathrm{SU}(3)$-homogeneous Aloff-Wallach's spaces. All these cases
belongs to (2) in Theorem \ref{Case-I-preparation-thm}.

W. Ziller pointed out that the fixed point set technique (i.e. the totally geodesic technique) in \cite{WZ2015} can be applied to discuss the unfinished (2) in
Theorem \ref{Case-I-preparation-thm}. Combining the fixed point set technique and the method from the homogeneous flag curvature formula, M. Xu
and W. Ziller proved following main theorem in \cite{XZi2016}.

\begin{theorem}\label{Case-I-thm}
Let $M=G/H$ be a compact simply connected homogeneous space
with a reversible Finsler metric with positive flag curvature on which
the compact Lie group $G$ acts by isometries. If $M$ is odd dimensional
and $\mathfrak{h}$ is a regular subalgebra of $\mathfrak{g}$, then either
$G/H$  carries a Riemannian homogeneous metrics with positive sectional
curvature, or $G/H$ is one of
\begin{description}
\item{\rm (1)} $\mathrm{Sp}(2)/\mathrm{diag}(z,z^3)$ with $z\in\mathbb{C}$;
\item{\rm (2)} $\mathrm{Sp}(2)/\mathrm{diag}(z,z)$ with $z\in\mathbb{C}$;
\item{\rm (3)} $\mathrm{Sp}(3)/\mathrm{diag}(z,z,q)$ with $z\in\mathbb{C}$, $q\in\mathbb{H}$;
\item{\rm (4)} $\mathrm{SU}(4)/\mathrm{diag}(zA,z,\bar{z}^3)$
with $A\in\mathrm{SU}(2)$ and $z\in\mathbb{C}$;
\item{\rm (5)} $\mathrm{G}_2/\mathrm{SU}(2)$ where $\mathrm{SU}(2)$ is the normal subgroup of $\mathrm{SO}(4)$ corresponding to the long root.
\end{description}
\end{theorem}

Finally, we remark that a general classification for odd dimensional
positively curved homogeneous Finsler spaces without the reversibility becomes
hopeful due to  the following theorem.

\begin{theorem}\label{powerful-criteria}
Let $F$ be a positively curved homogeneous Finsler metric on an odd dimensional coset space $G/H$. Then there does not exist a pair of linearly independent roots $\alpha$ and $\beta$ of $\mathfrak{g}$ such that the followings hold simultaneously:
\begin{description}
\item{\rm (1)} Neither $\alpha$ nor $\beta$ is a root of $\mathfrak{h}$.
\item{\rm (2)} $\pm\alpha$ are the only roots of $\mathfrak{g}$ in
$\mathbb{R}\alpha+\mathfrak{t}\cap\mathfrak{m}$.
\item{\rm (3)} $\pm\beta$ are the only roots of $\mathfrak{g}$
in $\mathbb{R}\alpha\pm\beta+\mathfrak{t}\cap\mathfrak{m}$.
\end{description}
\end{theorem}

Theorem \ref{powerful-criteria} is just  Lemma \ref{key-lemma-2} without  the reversible assumption. But the geometric phenomena behind these two statements are very different. In fact, in the proof of  them by deducing a contradiction, the methods to find   a zero flag curvature
 are essentially different.
Consequently the proof of Theorem \ref{powerful-criteria} is much longer and harder, in which we need  the more complicated homogeneous
flag curvature formula of L. Huang in \cite{Huang2013}. Theorem \ref{powerful-criteria} is
powerful for studying homogeneous positive flag curvature in that
many coset spaces can be excluded  from the list of possible candidates.

\subsection{Special classifications}
\label{s-c-t}

Now we turn to  the classification  of  non-Riemannian
homogeneous Randers or $(\alpha,\beta)$-spaces with positive flag curvature and vanishing
S-curvature, which can be summarized as the following theorem \cite{HD11} \cite{XD2015}.

\begin{theorem}\label{classification-1-with-0-S-curv}
Let  $G/H$ be a coset space which admits a
$G$-invariant non-Riemannian Randers or $(\alpha,\beta)$-metric with positive flag curvature
and vanishing S-curvature. Then  $G/H$ must be equivalent to
one of the following:
\begin{description}
\item{\rm (1)} The homogeneous spheres $S^{2n-1}=\mathrm{SU}(n)/\mathrm{SU}(n-1)$,
    $S^{2n-1}=\mathrm{U}(n)/\mathrm{U}(n-1)$,
    $S^{4n-1}=\mathrm{Sp}(n)/\mathrm{Sp}(n-1)$ and
    $S^{4n-1}=\mathrm{Sp}(n)\mathrm{U}(1)/\mathrm{Sp}(n-1)\mathrm{U}(1)$;
\item{\rm (2)} The Aloff-Wallach spaces
$S_{k,l}=\mathrm{SU}(3)/\mathrm{diag}(z^k,z^l,\bar{z}^{k+l})$ with $z\in\mathbb{C}$ or their $\mathrm{U}(3)$-homogeneous presentations.
\end{description}
\end{theorem}

The Randers case of Theorem \ref{classification-1-with-0-S-curv} was proved by Z. Hu and S. Deng in \cite{HD11}.
They studied non-Riemannian homogeneous Randers spaces with isotropic S-curvature (or equivalently vanishing S-curvature), and found that  these metrics
can be induced by a  {\it Killing navigation process}. To be precise, the non-Riemannian homogeneous Randers metric $F$ on $G/H$ corresponds to
the navigation datum $(F',V)$, where $F'$ is a Riemannian metric and $V$
is a vector field. Roughly speaking, it means each indicatrix $\{y\in T_xM|F(x,y)=1\}$ is a parallel shifting of the indicatrix
$\{y\in T_xM|F'(x,y)=1\}$ by the vector $V(x)$ (so we must have $F'(V(x))<1$ everywhere). In the homogeneous context,
$F'$ is a Riemannian homogeneous
metric on $G/H$, and $V$ is a $G$-invariant vector field on $G/H$. The vector field $V$
corresponds to a nonzero vector $v$ fixed by all $\mathrm{Ad}(H)$-actions.
Since the S-curvature vanishes,
the vector field $V$ involved in the navigation process is a Killing vector
field for $(G/H,F')$ as well as for $(G/H,F)$. That is the reason that we
call it a {\it Killing navigation process}.

For Killing navigation process, there is an important flag curvature equality given by the following theorem; see also \cite{HMO}.
\begin{theorem}\label{killing-navigation-flag-curvature-formula}
Let $F$ be the Finsler metric on $M$ induced by a navigation process with the navigation datum $(F',V)$, such that $V$ is a Killing vector field for $(M,F')$. Given any nonzero $y\in T_xM$, denote $\tilde{y}=y+{V(x)}/{F'(x,y)}$, i.e.,
$F'(x,y)=1$ iff $F(x,\tilde{y})=1$. Then for any linearly independent vectors
$y$ and $v$ satisfying $\langle y,v\rangle_y^{F'}=0$, we have the following
equality between flag curvatures,
$$K^{F'}(x,y,y\wedge v)=K^F(x,\tilde{y},\tilde{y}\wedge v).$$
\end{theorem}

Using Theorem \ref{killing-navigation-flag-curvature-formula}, we can show
that any non-Riemannian homogeneous Randers metric $F$ on $G/H$ with
positive flag curvature and vanishing S-curvature, corresponds to
a navigation datum $(F',V)$, such that $F'$ is positively curved Riemannian
homogeneous metric. By Theorem \ref{overall-classification-theorem-in-Riemannian-case}, only for those listed
in Theorem \ref{classification-1-with-0-S-curv} there exists nonzero vectors
in $\mathfrak{m}$ fixed by all $\mathrm{Ad}(H)$-actions. On the other hand,
all positively curved Riemannian homogeneous metrics on those coset spaces
have been explicitly described. So the homogeneous Randers metrics after the
Killing navigation process can be explicitly described as well.

 Now we consider the $(\alpha,\beta)$ case of Theorem \ref{classification-1-with-0-S-curv}.
For a non-Riemannian homogeneous $(\alpha,\beta)$-space
$(G/H,F)$ with vanishing S-curvature, we can find a closed subgroup $K$ in $G$, such that $\mathrm{Lie}(K)=\mathfrak{k}=\mathfrak{h}\oplus\mathbb{R}$,
and a Riemannian homogeneous metric $F'$ can be induced by $F$ and the Finsler
submersion $\pi:G/H\rightarrow G/K$. If $(G/H,F)$ is positively curved,
then so is $G/K$, which can be determined by N. Wallach's classification list \cite{Wallach1972}.

Finally, as an application of
Theorem \ref{Case-I-preparation-thm} and
Theorem \ref{Case-I-thm}, we consider the classification of non-Riemannian positively curved reversible homogeneous $(\alpha,\beta)$-spaces.  Notice that the proof of Theorem \ref{Case-I-thm}
uses the totally geodesic subspace technique to setup an induction. In each
step, the totally geodesic subspace remains to be a non-Riemannian positively curved reversible homogeneous $(\alpha,\beta)$-space when its dimension is bigger than $1$. So we only need to check the five undetermined cases in Theorem
\ref{Case-I-thm}. Three of them do not admit such an metric. It is interesting  that in dealing with
the coset space $\mathrm{Sp}(2)/\mathrm{diag}(z,z^3)$ with $z\in\mathbb{C}$,
 the argument in \cite{WZ2015} can be applied, which is an amazing application of Synge's theorem. The main results can be  summarized as the following theorem  \cite{XZh2016}.

\begin{theorem}\label{classification-reversible-alpha-beta}
Let $M=G/H$ be a compact simply connected homogeneous space which admits
a reversible $G$-invariant $(\alpha,\beta)$-Finsler metric with positive
flag curvature. Then either $G/H$  carries a Riemannian homogeneous metric
with positive sectional curvature, or it is one of the coset spaces
$\mathrm{Sp}(2)/\mathrm{diag}(z,z)$ or
$\mathrm{Sp}(3)/\mathrm{diag}(z,z,q)$ with $z\in\mathbb{C}$,
$q\in\mathbb{H}$.
\end{theorem}

\section{An explicit classification}\label{explicit}
As we have shown above, the classification of positively curved Finsler spaces is rather involved and difficult. However, there is a case in which we can give a complete classification of the metrics under isometries. We now recall the detail of this classification.

For simplicity, we will only present  the classification of non-Riemannian homogeneous  Randers metrics with isotropic S-curvature and positive flag curvature.
These metrics will be divided into five classes which can be described as the following.

\textbf{Case 1}\quad Invariant Randers metrics on the coset space $S^{2n+1}=\mathrm{SU}(n+1)/\mathrm{SU}(n)$. Denote ${\mathfrak g}=
\mathfrak{su}(n+1),{\mathfrak h}=\mathfrak{su}(n)$ and let
$\mathfrak{g}=\mathfrak{h}+\mathfrak{m}$ be the orthogonal decomposition of
$\mathfrak{g}$ with respect to the Killing form of $\mathfrak{g}$.
Then  ${\frak m}$
has a decomposition as
$$ {\mathfrak m}={\mathfrak m}_0\oplus{\mathfrak m}_1,$$
where
\begin{equation}\label{mz1}
    {\mathfrak m}_0={\mathbb R}X,\ \ X=\sqrt{-1}\left(
   \begin{array}{cc} -\frac{1}{n}E & 0 \\ 0 & 1 \end{array}
   \right),
\end{equation}
and
\begin{equation}
  {\mathfrak m}_1=\left\{\left.\left(
   \begin{array}{cc} 0 & \alpha \\ -\overline{\alpha}' & 0
   \end{array}
   \right)\right| \alpha '=(x_1,\cdots,x_n)\in {\Bbb C}^n\right\}.
  \end{equation}
An $\mathrm{SU}(n+1)$-invariant Riemannian metric on $\mathrm{SU}(n+1)/\mathrm{SU}(n)$ can be expressed as
\begin{equation}\label{mz2}
 h_t( X_1,X_2)
  =tc_1c_2+\mathrm{Re} (\alpha_1'\overline{\alpha_2}),
  \end{equation}
 with $t>0$, where $$X_i=c_i\sqrt{-1}\left(
   \begin{array}{cc} -\frac{1}{n}E & 0 \\ 0 & 1 \end{array}
   \right)+\left(
   \begin{array}{cc} 0 & \alpha_i \\ -\overline{\alpha}'_i & 0
   \end{array}
   \right),i=1,2.$$
This metric has positive sectional curvature if and only if
$0<t<\frac{2(n+1)}{3n}$.

Up to a scalar, there exists a unique  non-zero  $SU(n)$-invariant vector $X$ in  ${\frak m}$.  Then on the coset space $S^{2n+1}=SU(n+1)/SU(n)$ any $SU(n+1)$-invariant
Randers metric with almost isotropic S-curvature and positive flag curvature must be
a Randers metric which solves the Zermelo navigation problem of the Riemannian metric $h_t$
in (\ref{mz2}), with $0<t<\frac{2(n+1)}{3n}$, under the influence of the vector field generated
by the vector $cX$, where $X$ is defined in (\ref{mz1}),  and $|c|<\frac{1}{\sqrt{t}}$. Denote by the above metric by $F_{(t, c)}$. Then $F_{(t, c)}$ is non-Riemannian if and only if $c\ne 0$, and $F_{(t_1, c_1)}$ is isometric to
$F_{(t_2, c_2)}$ if and only if $t_1=t_2$ and $|c_1|=|c_2|$.

\textbf{Case 2}\quad Invariant Randers metrics on the coset space $S^{4n+3}=\mathrm{Sp}(n+1)/\mathrm{Sp}(n)$. The subspace ${\mathfrak m}$ of $\mathfrak{sp}(n+1)$ to be
$${\frak m}={\frak m}_0\oplus{\frak m}_1,$$
where
\begin{equation}\label{mz3}
{\mathfrak m}_0={\mathbb R}X_1\oplus {\mathbb R}X_2\oplus {\mathbb R}X_3
\end{equation}
is the subspace of $H$-fixed vectors in ${\mathfrak m}$, and
$X_i,i=1,2,3$, denote the elements of ${\mathbb H}^{(n+1)\times (n+1)}$
with the only non-zero element at the $(n+1,n+1)$-entry and  equal to
$\sqrt{2}I,\sqrt{2}J$ and $\sqrt{2}K$ respectively, here $I, J, K$ denote the standard imaginary units in ${\mathbb H}$, and
$$     {\mathfrak m}_1=\left\{\left.\left(
   \begin{array}{cc} 0 & \alpha \\ -\overline{\alpha}' & 0
   \end{array}
   \right)\right| \alpha '=(x_1,\cdots,x_n)\in {\mathbb H}^n\right\}.$$
Then, up to a positive multiple, any $\mathrm{Sp}(n)$-invariant Riemannian metric on $S^{4n+3}=\mathrm{Sp}(n+1)/\mathrm{Sp}(n)$ can be written as
\begin{equation}\label{mz4}
g_{(t_1, t_2, t_3)}(Y,Z)=t_1y_1z_1+t_2y_2z_2+t_3y_3z_3+{\rm Re}(\xi^*\eta),
\end{equation}
where $t_1, t_2, t_3$ are positive real numbers  and
\begin{eqnarray*}
Y &= & y_1X_1+y_2X_2+y_3X_3+\left(
   \begin{array}{cc} 0 & \xi \\ -\overline{\xi}' & 0
   \end{array}
   \right),\\
Z &=& z_1X_1+z_2X_2+z_3X_3+\left(
   \begin{array}{cc} 0 & \eta \\ -\overline{\eta}' & 0
   \end{array}
   \right).
\end{eqnarray*}
The condition for such a metric to have positive curvature can be stated as follows.
Let
$$V_i = (t^2_j+ t^2_k -3t^2_i+ 2t_it_j + 2t_it_k- 2t_jt_k)/t_i \,\,\mbox{ and}\,\, H_i = 4- 3t_i,$$
with $(i, j, k)$ a cyclic permutation of $(1, 2, 3)$. Then it is shown in \cite{VZ2009} that
 the homogeneous metrics $g_{(t_1,t_2,t_3)}$ on $S^{4n+3}$ have positive sectional curvature
if and only if
\begin{equation}\label{mz5}
V_i > 0 , H_i > 0\,\,\mbox{\rm and}\,\, 3|t_jt_k -t_j - t_k + t_i| < t_jt_k +\sqrt{H_iV_i},
\end{equation}
with $(i, j, k)$ a cyclic permutation of (1, 2, 3). It is also pointed out in \cite{VZ2009} that the set $(t_1, t_2, t_3)$ satisfying the above condition forms a non-empty slice.

  Now suppose $F$ is a Randers metric which solves the Zermelo's navigation problem of the
Rieamnnian metric $g_{(t_1, t_2, t_3)}$ in (\ref{mz4}), with $t_1\ne t_2=t_3\ne 1$,  under the influence of a vector field generated by $cX_1\in {\frak m}_0$, with $X_1$ as  in (\ref{mz3}),
and $0<|c|<\frac{1}{t_1}$. Then $F$ is a non-Riemannian invariant Randers metrics on $\mathrm{Sp}(n+1)/\mathrm{Sp}(n)$ with positive flag curvature and vanishing S-curvature. Two pairs $(g_{(t_1, t_2, t_3)}, cX_1)$ and $(g_{(t_1', t_2', t_3')}, c'X_1)$ correspond to isometric Randers metrics if and only if
$(t_1, t_2, t_3)=(t_1', t_2', t_3')$ and $|c|=|c'|$.

 \textbf{Case 3}\quad In this case, we also have $G/H=\mathrm{Sp}(n+1)/\mathrm{Sp}(n)$. Suppose   $F$ is a Randers metric which solves the Zermelo's navigation problem of the
Rieamnnian metric $g_{(t_1, t_2, t_3)}$ in (\ref{mz4}), with $t_1=t_2=t_3=t>0$,  under the influence of a vector field generated by $aX_1+bX_2+cX_3\in {\frak m}_0$ in (\ref{mz3}),
and $0<|a^2+b^2+c^2|<\frac{1}{t}$. Then  $F$ is also a non-Riemannian invariant Randers metrics on $\mathrm{SP}(n+1)/\mathrm{Sp}(n)$ with positive flag curvature and vanishing S-curvature. Two pairs $(g_{(t, t, t)}, aX_1+bX_2+cX_3)$ and $(g_{(t', t', t')}, a'X_1+b'X_2+c'X_3)$ correspond to isometric Randers metrics if and only if
$t=t'$ and $a^2+b^2+c^2=(a')^2+(b')^2+(c')^2$.

\textbf{Case 4}\quad  Consider the Aloff-Wallach's spaces $N_{(k,l)}=\mathrm{SU}(3)/S^1_{(k,l,-k-l)}$, with $kl(k+l)\ne 0$ and $\mathrm{gcd} (k, l)=1$£¬ where $S^1_{(k,l,-k-l)}=\{\mathrm{diag}(e^{k\theta\sqrt{-1}}, e^{l\theta\sqrt{-1}}, e^{-(k+l)\theta\sqrt{-1}})|\theta\in\mathbb{R}\}$.
We first consider the case $k\ne l$. In this case the isotropy representation has a decomposition as
\begin{eqnarray*}
{\mathfrak m}=V_0\oplus V_{2k+l}\oplus V_{2l+k}\oplus
V_{k-l},
\end{eqnarray*}
where
\begin{eqnarray*}
V_{2k+l} &=&\left.\left\{\left(
   \begin{array}{ccc} 0 & 0 & z\\ 0 & 0 & 0\\ -\overline{z} & 0 & 0 \end{array}
   \right)\right| z\in \mathbb{C}\right\},\\
V_{2l+k} &=& \left.\left\{\left(
   \begin{array}{ccc} 0 & 0 & 0\\ 0 & 0 & z\\ 0 & -\overline{z} & 0 \end{array}
   \right)\right| z\in \mathbb{C}\right\},\\
V_{k-l} &=&\left.\left\{\left(
   \begin{array}{ccc} 0 & z & 0\\-\overline{z} & 0 & 0\\ 0 & 0 & 0 \end{array}
   \right)\right| z\in \mathbb{C}\right\},
   \end{eqnarray*}
and
\begin{eqnarray*}
V_0=\left.\left\{\sqrt{-1}\left(
   \begin{array}{ccc} \ a & 0 & 0\\ 0 & b & 0\\ 0 & 0 & -(a+b) \end{array}
   \right)\right| (2k+l)a+(2l+k)b=0,a,b\in\mathbb{R}\right\}.
\end{eqnarray*}
The Lie algebra of $S_{k,l, -k-l}^1 =\{\mathrm{diag}(e^{k\theta\sqrt{-1}}, e^{l\theta\sqrt{-1}}, e^{-(k+l)\theta\sqrt{-1}})|\theta\in \mathbb{R}\}$ is ${\mathbb R}h_{k,l}$, where $h_{k,l}=\sqrt{-1}\mathrm{diag}(k,l,-(k+l))$. Set
$$V_1=V_0\oplus V_{k-l},V_2=V_{2k+l}\oplus
V_{2l+k},$$ and define
 $$q_0(X,Y)=-Re(tr(XY)), \quad X, Y\in \mathfrak{su}(3).$$
 Then  $q_0$ is an $\mathrm{Ad}(\mathrm{SU}(3)$-invariant inner product on $\mathfrak{su}(3)$, hence it  defines a   bi-invariant
 Riemannian metric on $\mathrm{SU}(3)$. It is easy to check that
$$q_0(h_{k,l},V_1)=0,\quad q_0(V_1,V_2)=0,$$
and
\begin{equation}\label{mz41}
[V_1,V_1]\subset {\mathbb R}h_{k,l}+V_1,[V_1,V_2]\subset V_2,[V_2,V_2]\subset
{\mathbb R}h_{k,l}+V_1.
\end{equation}
For $X, Y\in {\mathfrak m}$, set $X=X_1+X_2,Y=Y_1+Y_2,X_i,Y_i\in
V_i$, and define
\begin{equation}\label{mz42}
q_t(X,Y)=(1+t)q_0(X_1,Y_1)+q_0(X_2,Y_2)=q_0(X,Y)+tq_0(X_1,Y_1).
\end{equation}
It is shown in \cite{AW75} that if $-1<t<0$, then $q_t$ defines a $\mathrm{SU}(3)$-invariant Riemannian metric on
$N_{(k,l)}=\mathrm{SU}(3)/S^{1}_{(k,l, -k-l)}$ with positive curvature.

Now an $\mathrm{SU(3)}$-invariant non-Riemannian Randers metric $F$  on $N_{(k,l)}$ with $\gcd ( k,l)=1$, $kl(k+l)\ne 0$ and $(k,l)\ne (1,1)$ has positive flag curvature and vanishing S-curvature if and only if  $F$ is a Randers metric which solves the Zermelo's navigation problem of the
Rieamnnian metric $q_t$ in (\ref{mz42}), with $-1<t<0$,  under the influence of a vector field generated by an element
$w_a=\mathrm{diag} (a\sqrt{-1}, b\sqrt{-1}, -(a+b)\sqrt{-1})\in {\frak m}$, where $a, b$ are real numbers satisfying the conditions  $(2k+l)a+(2l+k)b=0$
and $ 0<a^2+b^2+(a+b)^2<\frac{1}{1+t}$. Two pairs $(h_t, w_a)$ and $(h_{t'}, w_{a'})$ correspond to isometric Randers metrics if and only if
$t=t'$ and $|a|=|a'|$.

\textbf{Case 5}\quad Consider the coset space $G/H=N_{(1,1)}=\mathrm{SU}(3)/S^1_{(1,1, -2)}$. The situation here is similar to the above case.  In this case, we set
$$V_0=\left.\left\{\sqrt{-1}\left(
   \begin{array}{ccc} \ a & 0 & 0\\ 0 & -a & 0\\ 0 & 0 & 0\end{array}
   \right)\right| a\in{\mathbb R}\right\},$$
and
\begin{eqnarray*}
V' &=&\left.\left\{\left(
   \begin{array}{ccc} 0 & 0 & z\\ 0 & 0 & 0\\ -\overline{z} & 0 & 0 \end{array}
   \right)\right| z\in \mathbb{C}\right\},\\
V'' &=& \left.\left\{\left(
   \begin{array}{ccc} 0 & 0 & 0\\ 0 & 0 & z\\ 0 & -\overline{z} & 0 \end{array}
   \right)\right| z\in \mathbb{C}\right\},\\
V''' &=&\left.\left\{\left(
   \begin{array}{ccc} 0 & z & 0\\-\overline{z} & 0 & 0\\ 0 & 0 & 0 \end{array}
   \right)\right| z\in \mathbb{C}\right\}.
   \end{eqnarray*}
Then the tangent space decomposes into the direct sum of irreducible
$S^{1}_{(1,1,-2)}$-submodules as
$${\mathfrak m}=V_0\oplus  V'\oplus V''\oplus V'''.$$
Define
$$V_1=V_0\oplus V''',\quad V_2=V'\oplus V''.$$
And define $q_0$, $q_t$ as in Case (4), then $q_t$ defines a $\mathrm{SU}(3)$-invariant Riemannian metric on
$N_{(1,1)}$, and this metric has positive curvature if and only if $-1<t<0$.

An $ \mathrm{SU}(3)$-invariant non-Riemannian Randers metric on $N_{(1,1)}$ has positive flag curvature and vanishing S-curvature if and only if  $F$ is a Randers metric which solves Zermelo's navigation problem of the Riemannian metric defined by $q_t$ in (\ref{mz42}), with
$t\in (-1,0)$, under the influence of a vector field generated by an element
$$w_{(a,b,c)}=\left(\begin{matrix}a\sqrt{-1} & b+c\sqrt{-1} & 0\cr
-b+c\sqrt{-1} & -a\sqrt{-1} & 0\cr
0 & 0 & 0\end{matrix}\right)\in {\mathfrak m},$$
where $a,b,c$ are real numbers satisfying the condition
$0<a^2+b^2+c^2 <\frac{1}{2(1+t)}$. Two pairs $(h_t, w_{(a,b,c)})$ and $(h_{t'}, w_{(a', b', c')})$ correspond to isometric Randers metrics  if and only if $t=t'$ and
$a^2+b^2+c^2=(a')^2+(b')^2+(c')^2$.

We summarize the above as the following theorem.
\begin{theorem} (\cite{HD11})\quad
Let $(G/H, F)$ be a homogeneous non-Riemannian Randers metric. Then $F$ has positive flag curvature and vanishing S-curvature if and only if it is one of
the metrics in the above cases (1)-(5). Two Randers metrics in different cases in the above can not  be isometric.
\end{theorem}

\section{The flag-wise positively curved condition}

In Finsler geometry, we may study other positive curvature conditions. Due to the  dependence on the
flag pole,  flag curvature is much more local than
sectional curvature in Riemannian geometry. This motivates us to define  the following {\it flag-wise positively curved condition}, or {\it (FP) Condition} for short \cite{XD-FP}.

\begin{definition}
A Finsler space $(M,F)$ is called flag-wise positively curved, or satisfying the (FP) Condition, if for any $x\in M$ and any tangent plane $\mathbf{P}\subset
T_xM$, there exists  a nonzero tangent vector $y\in\mathbf{P}$ such that
$K^F(x,y,\mathbf{P})>0$.
\end{definition}

Notice that in Riemannian geometry, (FP) Condition  is equivalent to the positively curved condition. So
among all non-negatively curved metrics, those satisfying (FP) Condition
provide a very special class, which only appears in Finsler geometry.
Strictly speaking, the combination of non-negatively curved condition with (FP) Condition is still weaker than the positive flag curvature condition. However,
intuitively it is already  "sufficiently" close to the positive flag curvature condition.  This idea leads to a number of explicit examples showing the differences, and some consideration on
 the geometric properties of compact or complete Finsler spaces satisfying flag-wise positively curved and non-negatively curved conditions at the same time.

In \cite{XD-FP}, we studied these problems in homogeneous Finsler geometry and proved
the following theorem.

\begin{theorem} \label{FP-theorem}
There are many compact coset spaces $G/H$  which admit
flag-wise positively and non-negatively curved $G$-homogeneous Finsler metrics,
but cannot be endowed with any positively curved $G$-homogeneous Finsler metric.
\end{theorem}

The $S^1$-bundles over compact Hermitian symmetric spaces provides many coset spaces possessing the properties in  Theorem \ref{FP-theorem}. Using Theorem \ref{powerful-criteria}, we can prove that most of these coset spaces
do not admit positively curved homogeneous Finsler metrics.

The examples considered in \cite{XD-FP} are canonical $S^1$-bundles over compact irreducible
Hermitian symmetric spaces, i.e.,
\begin{eqnarray}\label{FP-list-1}
& &\mathrm{SU}(p+q)/\mathrm{SU}(p)\mathrm{SU}(q),\, \mathrm{SO}(n)/\mathrm{SO}(n-2),\,
\mathrm{Sp}(n)/\mathrm{SU}(n),
\nonumber\\
& &\mathrm{SO}(2n)/\mathrm{SU}(n),\,
\mathrm{E}_6/\mathrm{SO}(10)\mbox{ and }\mathrm{E}_7/\mathrm{E}_6.
\end{eqnarray}

The corresponding metric $F$ in each case is a homogeneous Randers metric corresponding to a navigation datum $(F',V)$, where $F'$ is a Riemannian
normal homogeneous metric, and $F'$ is an invariant Killing vector field.
By Theorem \ref{killing-navigation-flag-curvature-formula}, $F$ is non-negatively curved. On the other hand, the flags and poles with zero flag curvature are rearranged by the Killing navigation process, which make a chance for the (FP) Condition.

Theorem \ref{powerful-criteria} can be used  to show that most coset spaces in
(\ref{FP-list-1}) do not admit positive flag curvature. To be precise, the
following coset spaces do not admit homogeneous Finsler metrics with positive flag curvature,
\begin{eqnarray*}
& &\mathrm{SU}(p+q)/\mathrm{SU}(p)\mathrm{SU}(q), \mbox{ with }p>q\geq 2
\mbox{ or }p=q>3,\nonumber\\
& &\mathrm{Sp}(n)/\mathrm{SU}(n),\mbox{ with }n>4,\nonumber\\
& &\mathrm{SO}(2n)/\mathrm{SU}(n), \mbox{ with }n=5\mbox{ or }n>6,
\nonumber\\
& &\mathrm{E}_6/\mathrm{SO}(10)\mbox{ and }\mathrm{E}_7/\mathrm{E}_6,
\end{eqnarray*}
which meet the requirement of Theorem \ref{FP-theorem}.

The classification of flag-wise positively and non-negatively curved homogeneous Finsler spaces seems to be an interesting project. We conjecture that
 any such a coset space  must be compact and  satisfies the rank inequality (see (2) of Theorem \ref{rank-inequ-thm}).

Finally,  we remark that the (FP) Condition  is a weak condition. In fact, in many cases, a  perturbation of
a non-negatively curved Finsler metric  will produce
a flag-wise positively curved Finsler metric. More  precisely, we have
the following theorem  \cite{X}.

\begin{theorem}\label{thm-5-3}
Let $G/H$ be a compact coset space with a finite fundamental group. Then
$G/H$ and $G/H\times S^1$ admits (generally non-homogeneous) Finsler metrics
satisfying the (FP) Condition.
\end{theorem}

The metrics predicted by Theorem \ref{thm-5-3} are constructed by a
combination of the Killing navigation process and a gluing technique.
By similar constructive method, we classified all
quasi-compact Lie groups admitting a flag-wise positively curved left invariant Finsler metrics, namely, we have the following theorem \cite{X}.

\begin{theorem}
A connected quasi-compact Lie group $G$ admits flag-wise
positively curved left invariant Finsler metrics if and only if  $\dim C(G)<2$.
\end{theorem}

Notice that if  $\dim C(G)\geq 2$, then the homogeneous flag curvature
formula (\ref{homogeneous-flag-curvature-formula}) implies that, for any left invariant Finsler
metric $F$, any tangent plane $\mathbf{P}\subset\mathfrak{c}(\mathfrak{g})\subset\mathfrak{g}=T_eG$
must have zero flag curvature, no matter which flag pole we choose.

As an interesting application of the (FP) Condition, we prove the following theorem:
\begin{theorem}(\cite{XD-FP})\quad On the product manifolds $S^2\times S^3$ and $S^6\times S^7$, there exist non-negatively curved Finsler metrics which satisfy the (FP) Condition.
\end{theorem}

It is a well known long standing open problem whether there exists a product manifold which admits a positively curved Riemannian metrics, which is a part of the generalized Hopf's conjecture. The interest of the above theorems lies in the fact that at least there exist some product manifolds
which admit Finsler metrics of non-negative flag curvature, and for any tangent plane there exists a flag pole in that plane such that the corresponding flag curvature is positive.

\section{Homogeneous Finsler spaces of negative curvature}
As a counterpart of the study of homogeneous Finsler spaces of positive flag curvature, the study of the negatively curved ones has also led to a series of interesting results. Contrary to the positive case, in the negative we got a lot of rigidity results. We now give a survey of the main progress.

The following theorem was proved in \cite{DHneg1}
\begin{theorem}\label{nega1}
Let $(M, F)$ be a connected
homogeneous Finsler space of non-positive flag curvature. If the
Ricci scalar is everywhere strictly negative, then $M$ is simply
connected.
\end{theorem}

As a special case of the above theorem, a connected homogeneous Finsler space of negative curvature must be simply connected. Note that the Riemannian case of
Theorem \ref{nega1} is proved by
S. Kobayashi in \cite{K}.

In the following we give some rigidity results.
\begin{theorem}{\rm (\cite{DHneg})}\quad
Let $(M, F)$ be a homogeneous Randers space. If $F$ is an Einstein metric and Ricci scalar is negative, then $F$ must be Riemannian.
\end{theorem}
In particular, a homogeneous Einstein-Randers metric with negative curvature must be Riemannian. This result is a special case of the following
\begin{theorem}\label{neg_ricci}
{\rm (\cite{DHneg1})}\quad
Let $(M, F)$ be a connected homogeneous Randers space. If $F$ has almost isotropic S-curvature and negative Ricci scalar, then $F$ must be Riemannian.
\end{theorem}

\section{Some Open problems}
In this section we collect some problems related to the topics of this paper.
\begin{problem}
Complete the classification of coset spaces $G/H$ which admit $G$-invariant Finsler metrics of positive flag curvature.
\end{problem}
 As  stated in the previous sections, the even dimensional case has been completely settled. In the odd dimensional case, great progress has been made by Xu-Deng \cite{XD2016} and Xu-Ziller \cite{XZi2016}. However, there are some special cases in which all the techniques available are not effective, and these cases are   the most difficult part of the problem. Note that in the general sense, we can only give a classification of the coset spaces, but it is impossible to give a complete classification of all the metrics under isometries. In fact, this can be achieved only in some very special cases, as explained in Section 8.

 \begin{problem}
 Classify homogeneous Randers spaces with positive flag curvature under isometries.
 \end{problem}

 In Section \ref{explicit}, we describe the main results of Hu-Deng \cite{HD11} on the classification of homogeneous Randers spaces with positive flag curvature and vanishing S-curvature. It is interesting to generalize the classification
 to the general Randers spaces without the S-curvature restriction. Note that an explicit coordinate-free formula of the flag curvature of homogeneous Randers spaces is given by Deng-Hu in \cite{DH2013}, and this formula will be useful in the study of this problem.

 \begin{problem}
 Classify homogeneous Finsler spaces with constant flag  curvature.
 \end{problem}
 The study of Finsler spaces with constant flag curvature has been one of the central problems in this field. The Randers case has been completely settled by Bao-Robles-Shen in \cite{BRS04}, and the homogeneous ones among their classification were figured out by Deng in \cite{De2010}. However, a general classification seems to be unreachable. It is easily seen that a homogeneous Finsler space with negative constant flag curvature must be Riemannian, and a homogeneous flat Finsler space must be locally Minkowskian (see, for example \cite{DE12}). Therefore one only need to study this problem for positive constant case. As pointed out  in \cite{CS2005}, a connected simply connected Finsler space with positive constant flag curvature must be differmorphic to a  sphere.
 Therefore a connected simply connected homogeneous Finsler space with positive constant flag curvature can be viewed as an invariant Finsler metric on a coset space $G/H$, where $G$ is a connected compact Lie group which has an effective transitive action on a sphere. Note that the connected compact Lie groups which admit an effective and transitive action have been classified by Montgomery-Samelson \cite{MS43} and Borel \cite{Bo1940}, and a list can be found in \cite{HD11}. It is easy to check that in almost all the cases the invariant Finsler metrics are Randers metrics. The difficult part of the problem is probably the classification of $\mathrm{Sp}(n+1)$-invariant Finsler metrics with positive constant flag curvature on the
 spheres $S^{4n+3}=\mathrm{Sp}(n+1)/\mathrm{Sp}(n)$.

 \begin{problem}
 Does there exist examples of non-Riemannian homogeneous Finsler spaces with negative Ricci scalar and vanishing S-curvature? If so, classify them.
 \end{problem}
By Theorem \ref{neg_ricci}, in the Randers case, the answer to the above problem is negative. If we drop the assumption of S-curvature, then there exists many examples in the Randers case. In fact, any homogeneous Riemannian manifold with negative curvature can be realized as a solvable Lie group endowed with a left invariant metric. Therefore, if one fix a tangent vector at the unit element with length sufficiently small, then using the navigation method, one can produce a left invariant Randers metric with negative flag curvature. However, the S-curvature of such a metric cannot be vanishing everywhere. On the other hand, if $F$
is a homogeneous Finsler space with negative flag curvature and vanishing S-curvature, then it must be Riemannian; this is a corollary of  a result of Shen \cite{SH1}; see \cite{DE12}. We conjecture that the answer to this problem is negative, namely, any homogeneous Finsler space with negative Ricci curvature and vanishing S-curvature must be Riemnanian.

\end{document}